\documentclass[12pt]{amsart}
\usepackage{amssymb}

\usepackage{mathptmx}

\newcommand{\mm}{\mathfrak m}

\newcommand{\Z}{\mathbb{Z}}

\DeclareMathOperator{\pnt}{\raise 0.5mm \hbox{\large\bf.}}
\DeclareMathOperator{\lpnt}{\hbox{\large\bf.}}

\DeclareMathOperator{\depth}{depth}

\DeclareMathOperator{\Ker}{Ker}

\DeclareMathOperator{\reg}{reg}
\DeclareMathOperator{\Tor}{Tor}

\def\+#1{\relax\ifmmode\if\noexpand #1\relax \mathop{\kern
    0pt^+{#1}}\nolimits\else \kern 0pt^+\!#1 \fi\else$^*$#1\fi}

\DeclareMathOperator{\Hom}{Hom}

\DeclareMathOperator{\pd}{pd}

\DeclareMathOperator{\Image}{Im}

\newtheorem{thm}{Theorem}[section]
\newtheorem{lem}[thm]{Lemma}

\newtheorem{quest}[thm]{Question}

\theoremstyle{definition}
\newtheorem{defn}[thm]{Definition}

\newtheorem{rem}[thm]{Remark}
\newtheorem{ex}[thm]{Example}

\theoremstyle{plain}
\newtheorem*{thm*}{Theorem}

\let\phi=\varphi
\title{On the regularity over positively graded algebras}

\author{Tim R\"omer}
\address{FB Mathematik/Informatik, Universit\"at Osnabr\"uck, 49069 Osnabr\"uck, Germany}
\email{troemer@uos.de}

\begin{document}

\begin{abstract}
We study the relationship between the Tor-regularity and the local-regularity
over a positively graded algebra defined over a field
which coincide if the algebra is a standard graded polynomial ring.
In this case both are characterizations of the so-called Castelnuovo--Mumford regularity.
Moreover,  we can characterize a standard graded
polynomial ring as a $K$-algebra with extremal properties with respect
to the
Tor- and the local-regularity.
For modules of finite projective dimension we get a nice formula
relating the two regularity notions.
Interesting examples are given to help to understand the relationship between
the Tor- and the local-regularity in general.
\end{abstract}

\maketitle
%
%
%
\section{Introduction}
Let $K$ be a field and $S=K[x_1,\dots,x_n]$ be a standard graded polynomial ring
with unique graded maximal ideal $\mm=(x_1,\dots,x_n)$.
Eisenbud and Goto \cite{EIGO} proved that
for a finitely generated graded $S$-module $M$
the finite numbers
\begin{enumerate}
\item
$\inf \{ r \in \Z : \text{ for all } i\geq 0 \text{ and all } s>r
\text{ we have } \Tor_{i}^{S}(M,K)_{i+s}=0 \}$
\item
$\inf\{ r\in \Z : \text{ for all } i\geq0 \text{ and all }s>r
\text{ we have }H_{\mm}^{i}(M)_{s-i}=0\}$
\end{enumerate}
coincide. Usually one calls this number the {\em Castelnuovo--Mumford}-regularity $\reg_S(M)$
of $M$.
Moreover, if
we denote by $M_{\geq q}$ for $q \geq 0$ the truncation  of $M$ defined as the graded $S$-module
with homogeneous components:
$(M_{\geq q})_{i}=M_{i}$ if $i\geq q$ and
$(M_{\geq q})_{i}=0$ for $i<q$, then
$\reg_S(M)$ is also the least number $q$ such that $M_{\geq q}$
is non trivial and has a $q$-linear $S$-resolution, i.e. $\dim_K \Tor_i^S(M_{\geq q},K)_{i+j}=0$
for $j\neq q$.

It is a natural question to understand the relationship
between these numbers in the situation where the $K$-algebra is not longer a polynomial ring.
In the following a
{\em positively graded $K$-algebra} $R$ is
a Noetherian commutative $K$-algebra
such that $R=\bigoplus_{i\geq 0} R_i$ with $R_0=K$.
We denote always by $\mm=\bigoplus_{i> 0} R_i$ the unique graded maximal ideal of $R$.
We say that $R$ is {\em standard graded} if $R$ is generated in degree $1$.
If $R$ is a polynomial ring, then we call $R$ a positively graded polynomial ring and standard graded polynomial ring respectively.
A finitely generated graded $R$-module is always a
non-trivial $\Z$-graded $R$-module $M=\bigoplus_{i\in \Z} M_i$.
The crucial definitions of this paper are:

\begin{defn}
\label{regt}
Let $R$ be a positively graded $K$-algebra and $M$ be a finitely generated graded $R$-module.
Then:
\begin{enumerate}
\item
$
\reg_{R}^{T}(M)
=
\inf \{ r \in \Z : \text{ for all } i\geq 0 \text{ and all } s>r
\text{ we have } \Tor_{i}^{R}(M,K)_{i+s}=0 \}
$
is called the
{\em Tor-regularity} of $M$.
\item
$
\reg_{R}^{L}(M)=
\inf\{ r\in \Z : \text{ for all } i\geq0 \text{ and all }s>r
\text{ we have }H_{\mm}^{i}(M)_{s-i}=0\}
$
is called the {\em local-regularity} of $M$.
\end{enumerate}
\end{defn}
After some preliminary remarks in Section \ref{prelim}
we prove in Section \ref{comparethenumbers} our
first main result:
\begin{thm}
Let $R$ be a positively graded $K$-algebra and $M$ be a finitely generated graded $R$-module.
Then:
\begin{enumerate}
\item
$\reg_{R}^{L}(M)-\reg_{R}^{L} (R) \leq \reg_{R}^{T} (M)$.
\item
If $R$ is standard graded, then
$\reg_{R}^{T} (M) \leq \reg_{R}^{L} (M) + \reg_{R}^{T} (K)$.
\end{enumerate}
\end{thm}
Observe that the upper inequality is essentially due to Avramov and Eisenbud \cite{AVEI92}.
J{\o}rgensen \cite{JO} proved a much more generally version for complexes over not necessarily commutative $K$-algebras
which have a balanced dualizing complex.
For modules we present here a straight forward proof
which avoids the technical machinery used in \cite{JO}.
See Herzog and Restuccia \cite{HERE} for a similar result over
standard graded $K$-algebras.
Note that by a result of Avramov and Eisenbud \cite{AVEI92}
if $R$ is a Koszul algebra, i.e. $\reg_{R}^{T}(K)=0$
where we consider $K=R/\mm$ naturally as an $R$-module,
it is still true that $\reg_{R}^{T}(M)$ is
the least number $q$ such that $M_{\geq q}$ is non trivial and has a $q$-linear $R$-resolution.

Having certain inequalities of invariants related to a module $M$,
it is of course interesting to understand for which modules equality holds.
Considering standard graded $K$-algebras $R$ we know by the graded version of the famous result
of Auslander--Buchsbaum--Serre that $R$ is a polynomial ring if and only if
$\pd_R (M)<\infty$ for all finitely generated graded $R$-modules.
Moreover, it is enough to show that $\pd_R (K)<\infty$ to conclude that $R$ is a polynomial ring.
Interestingly a polynomial ring is also characterized by extremal properties
with respect to the regularity notions introduced above.
More precisely, in Section \ref{borderline} we prove:
\begin{thm}
Let $R$ be a standard graded $K$-algebra. The following statements are equivalent:
\begin{enumerate}
\item
For all finitely generated graded $R$-modules $M$ we have $\reg^L_R (M)-\reg^L_R (R)=\reg^T_R (M)$;
\item
For all finitely generated graded $R$-modules $M$ we have $\reg^T_R (M)=\reg^L_R (M)+\reg^T_R (K)$;
\item
For all finitely generated graded $R$-modules $M$ we have $\reg^T_R (M)=\reg^L_R (M)$;
\item
$R$ is Koszul and $\reg^L_R (R)=0$;
\item
$R=K[x_1,\dots,x_n]$ is a standard graded polynomial ring.
\end{enumerate}
\end{thm}

In the general case we can still show the following nice fact:

\begin{thm}
Let $R$ be a positively graded $K$-algebra and $M$ be a finitely generated graded $R$-module
such that $\pd_R(M)<\infty$.
Then
$$\reg_{R}^{L}(M) - \reg_{R}^{L} (R) = \reg_{R}^{T}(M).$$
\end{thm}
Also see Chardin \cite{CH07} for similar results.
By giving an example that the converse of the latter result does not hold,
it still interesting to understand for which modules
we have
$\reg^L_R (M)-\reg^L_R (R)=\reg^T_R (M)$ and
$\reg^T_R (M)=\reg^L_R (M)+\reg^T_R (K)$ respectively.
We conclude the paper in Section \ref{exampleseries}
with the observation
that there  exists a Koszul algebra $R$ such that $\depth(R)>0$ and $r=\reg_{R}^{L}(R)>0$
and  we have for $0<j<r$ that
$$
0=\reg_{R}^{L}(\mm^j)-r < \reg_{R}^{T}(\mm^j)=j <r=\reg_{R}^{L}(\mm^j).
$$
In this sense any number between
$\reg^L_R (M)-\reg^L_R (R)$ and
$\reg^L_R (M)+\reg^T_R (K)$ can be the Tor-regularity of a module.

We are grateful to Prof.\ J. Herzog
for inspiring discussions on the subject of this paper.

%
%
%
\section{Preliminaries}
\label{prelim}
In this section we fix some further notation and recall some definitions.
For facts related to commutative algebra
we refer to the book of Eisenbud \cite{EI}.
A standard reference on homological algebra is Weibel \cite{WE95}.
Now following Priddy \cite{PR} we define:
\begin{defn}
\label{koszul}
Let $R$ be a standard graded $K$-algebra. Then $R$ is called a {\em Koszul algebra}
if $\reg_{R}^{T}(K)=0$ where we consider $K=R/\mm$ naturally as an $R$-module.
\end{defn}
For a positively graded $K$-algebra $R$
and a finitely generated graded $R$-module $M$ we say that $M$ has an
{\em $q$-linear resolution} if $\Tor_i^R(M,K)_{i+j}=0$ for $j\neq q$.
Thus if
we consider the minimal graded free resolution
$$
F_{\lpnt}:
\ldots\to F_i \overset{\partial_{i}}{\to} \ldots \overset{\partial_{1}}{\to} F_0 \to
M \to 0
$$
of $M$ with
$F_i=\bigoplus_{ j \in \Z } R(-j)^{\beta_{ij}^R(M)}$
where $\beta_{ij}^R(M) = \dim_K \Tor_l^R(M,K)_j$ are the {\em graded Betti-numbers} of $M$,
then $M$ has a $q$-linear resolution if and only if $\beta_{ii+j}^R(M)=0$ for $j \neq q$.
In particular, if $R$ is standard graded, then $R$ is Koszul if and only if
$K$ has a $0$-linear resolution. E.g.\ a standard graded polynomial ring $K[x_1,\dots,x_n]$ is trivially a Koszul algebra
since the Koszul complex on the variables $x_1,\dots,x_n$
provides a minimal graded free resolution of $K$ which is $0$-linear.
We will need the following result which is essentially due to
Avramov and Eisenbud \cite{AVEI92}:
\begin{thm}
\label{koszulavei}
Let $R$ be a Koszul algebra and $M$ be a finitely generated graded $R$-module, then
$$
\reg_R^T (M) \leq \reg_R^L (M)<\infty.
$$
\end{thm}
\begin{proof}
Let $R=S/I$ where $S$ is a standard graded polynomial ring
and $I \subset S$ is a graded ideal containing no linear forms.
Avramov and Eisenbud proved that
$\reg_R^T (M) \leq \reg_S^T (M)$.
But over a polynomial ring we have
$\reg_S^T (M) =\reg_S^L (M)$ by Eisenbud and Goto \cite{EIGO}.
Moreover, it is well-known that
$\reg_S^L (M)=\reg_R^L (M)$ simply because the local cohomology of $M$
with respect to the maximal ideal computed over $S$ is isomorphic to
the local cohomology of $M$
with respect to the maximal ideal computed over $R$.
That $\reg_R^L (M)<\infty$ follows now from the fact that
$\reg_S^L (M)<\infty$.
(E.g.\ see \cite{BRHE98}:
There are only finitely many  local cohomology groups not zero and
all of them have the property that $H^i_\mm(M)_j=0$ for $j \gg 0$.)
\end{proof}

The Koszul property can not be decided by knowing only finitely many
graded Betti numbers of $K$. (See \cite{RO} for examples.)
Recently Avramov and Peeva \cite{AVPE01}
proved the following remarkable result:
\begin{thm}
\label{koszulavpe}
Let $R$ be a positively graded $K$-algebra. Then the following statements are equivalent:
\begin{enumerate}
\item
$R$ is Koszul;
\item
$R$ is standard graded and for every finitely generated graded $R$-module $M$ we have
$\reg_R^T (M)<\infty$;
\item
$R$ is standard graded and we have $\reg_R^T (K)<\infty$.
\end{enumerate}
\end{thm}
Thus $K$ is a test-module for the Koszul property using the invariant $\reg_R^T (K)$.
In the next sections we will use occasionally the
following observations:

\begin{rem}
\label{goodtoknow}
Let $R$ be a positively graded $K$-algebra and $M$ be
a finitely generated graded $R$-module.
\begin{enumerate}
\item
We have $\reg^L_R (M)< \infty$. Indeed, the arguments used in the proof of
Theorem \ref{koszulavei} to show that $\reg^L_R (M)$ is finite,
did not used the fact that $R$ is Koszul.
\item
$\reg^L_R (K)=0$ because $H^i_\mm(K)=0$ for $i\neq 0$
and $H^0_\mm(K)=K$.
\item
By a result of Grothendieck (e.g. see \cite{BRHE98})
we know that
$$
H_{\mm}^{i}(M)
\begin{cases}
=0 & \text{ for }i<\depth (M) \text{ and } i>\dim (M), \\
\neq 0 & \text{ for }i=\depth (M) \text{ and }i= \dim (M).
\end{cases}
$$
\end{enumerate}
\end{rem}

%
%
%
\section{Comparison of the Tor- and the local-regularity}
\label{comparethenumbers}
We want to compare the notion of regularities as introduced in Section \ref{prelim}.
The main result of this section is the next theorem.

\begin{thm}
\label{haupt}
Let $R$ be a positively graded $K$-algebra and $M$ be a finitely generated graded $R$-module.
Then:
\begin{enumerate}
\item
$\reg_{R}^{L}(M)-\reg_{R}^{L} (R) \leq \reg_{R}^{T} (M)$.
\item
If $R$ is standard graded, then
$\reg_{R}^{T} (M) \leq \reg_{R}^{L} (M) + \reg_{R}^{T} (K)$.
\end{enumerate}
\end{thm}
\begin{proof}
(i):
If $\reg_{R}^{T} (M)=\infty$, then nothing is to show. Next assume that
$\reg_{R}^{T} (M)<\infty$.
Observe that the numbers $\reg_{R}^{L}(M)$ and $\reg_{R}^{L} (R)$ are always finite by Remark \ref{goodtoknow}.
We consider the minimal graded free resolution
$$
F_{\lpnt}:
\ldots\to F_l \overset{\partial_{l}}{\to} \ldots \overset{\partial_{1}}{\to} F_0 \to
M \to 0
$$
of $M$ with
$F_l=\bigoplus_{ j \in \Z } R(-j)^{\beta_{lj}^R(M)}$
where $\beta_{lj}^R(M) = \dim_K \Tor_l^R(M,K)_j$ are the graded Betti-numbers of $M$.
Note that $\beta_{lj}^R(M)=0$ for $j>l+\reg_{R}^{T} (M)$.
Define
$$
C_{l}:=\Ker \partial_{l} \text{ for } l\geq 0 \text{ and set } C_{-1}:=M.
$$
For $l\geq 0$ the short exact sequences
$$
0 \to C_{l} \to F_l \to C_{l-1} \to 0
$$
give rise to a long exact local cohomology sequence in degree $j-i$
$$
0
\to
H_{\mm}^{0}(C_{l})_{j-i}
\to
H_{\mm}^{0}(F_l)_{j-i}
\to
H_{\mm}^{0}(C_{l-1})_{j-i}
\to
\ldots
$$
$$
\to
H_{\mm}^{i}(C_{l})_{j-i}
\to
H_{\mm}^{i}(F_l)_{j-i}
\to
H_{\mm}^{i}(C_{l-1})_{j-i}
\to
\ldots
$$
Since $H_{\mm}^{i}(R(-j))=H_{\mm}^{i}(R)(-j)$ and
$H_{\mm}^{i}(\cdot)$ is an additive functor,
we have that for $j>\reg_{R}^{L} (R)+\reg_{R}^{T}(M) +l$
and for all $i\geq 0$ that
$$
H_{\mm}^{i}(F_l)_{j-i}
=
0.
$$
Thus for $l=0$ and $j>\reg_{R}^{L} (R)+\reg_{R}^{T}(M)$
we obtain
$$
H_{\mm}^{i}(M)_{j-i} \cong H_{\mm}^{i+1}(C_{0})_{j-i}.
$$
For $j+1>\reg_{R}^{L} (R)+\reg_{R}^{T}(M)+1 \Leftrightarrow j>\reg_{R}^{L} (R)+\reg_{R}^{T}(M)$
we get analogously
$$
H_{\mm}^{i+1}(C_{0})_{j-i} \cong H_{\mm}^{i+2}(C_{1})_{j-i}.
$$
Using an appropriate induction
we see that for $j>\reg_{R}^{L} (R)+\reg_{R}^{T}(M)$
we have
$$
H_{\mm}^{i}(M)_{j-i} \cong \ldots \cong H_{\mm}^{i+l+1}(C_{l})_{j-i}.
$$
Note that $\dim C_{l}\leq \dim R=:d$ and we get from Remark \ref{goodtoknow} that for $l\geq d-i$
we have $H_{\mm}^{i+l+1}(C_{l})=0$.
All in all we obtain for $j>\reg_{R}^{L} (R)+\reg_{R}^{T}(M)$ and $l\geq d-i$
that
$$
H_{\mm}^{i}(M)_{j-i}\cong H_{\mm}^{i+l+1}(C_{l})_{j-i}=0
$$
Hence $\reg_{R}^{L} (R)+\reg_{R}^{T}(M) \geq \reg_{R}^{L}(M) $ as desired.

(ii):
If $\reg_{R}^{T} (K)=\infty$, then nothing is to show.
Next assume $\reg_{R}^{T} (K)<\infty$.
Then it follows from Theorem \ref{koszulavpe} that
$\reg_{R}^{T} (K)=0$ and $R$ is a Koszul algebra.
But now the inequality
$\reg_{R}^{T} (M)
\leq
\reg_{R}^{L} (M)$ was shown in Theorem \ref{koszulavei}.
\end{proof}

In Section \ref{borderline} we will see that most times
$\reg_{R}^{L}(M) \neq \reg_{R}^{T} (M)$, so these two regularities do
no coincide in general.
For Koszul algebras we still have
the result that the regularity is related to linear resolutions of truncations of $M$.
(See \cite{AVEI92}, \cite{EIGO} and \cite{JO}.)
Here for a graded $R$-module $M$ and an integer $q$ we define
the truncation $M_{\geq q}$ of $M$ as the graded $R$-module
with homogeneous components:
$$
(M_{\geq q})_{i}=
\begin{cases}
M_{i} & \text{if $i\geq q$,} \\
0 & \text{else.}
\end{cases}
$$

\begin{thm}
\label{first}
Let $R$ be a Koszul algebra,
$M$ be a finitely generated graded $R$-module and $q \in \Z$.
The following statements are equivalent:
\begin{enumerate}
\item
$q\geq \reg_{R}^{T}(M)$;
\item
$\Tor_{i}^{R}(M,K)_{i+j}=0$ for all $i \geq 0$ and all $j>q$;
\item
$M_{\geq q}$ has a $q$-linear $R$-resolution.
\end{enumerate}
In particular,
$\reg_{R}^{T} (M)$ is the least  $q \in \Z$ such that $M_{\geq q}$ is non trivial and has a $q$-linear free resolution.
Moreover, if $q\geq \reg_{R}^{T} (M)$, then $M_{\geq q}$ has a $q$-linear free resolution.
\end{thm}
\begin{proof}
The equivalence of (i) and (ii) follows directly from the definition of $\reg_{R}^{T}(M)$.

Let now $F_{\lpnt}$ be a minimal graded free resolution of $K$ as an $R$-module.
Since $R$ is Koszul we have  $\reg_{R}^{T}(K)=0$ and thus
$
0 = \dim_{K} \Tor_{i}^{R}(K,K)_{i+j}=\beta^R_{i,i+j}(K)\text{ for } j\neq 0.$
Hence
$$F_{\lpnt}: \ldots \rightarrow R(-i)^{c_{i}} \rightarrow \ldots \rightarrow R^{c_{0}} \rightarrow K \rightarrow 0.$$
Assume (ii) holds.
The $K$-vector space
$\Tor_{i}^{R}(M,K)_{i+j}$ is the $i$-th homology of the following complex:
$$F_{\lpnt}\otimes_{R} M:
\ldots \rightarrow (R(-i)^{c_{i}}\otimes_{R} M)_{i+j} \rightarrow \ldots
\rightarrow (R^{c_{0}}\otimes_{R} M)_{i+j} \rightarrow 0.
$$
For $j>q$ we have
$
(R(-i)^{c_{i}}\otimes_{R} M)_{i+j}  =(R(-i)^{c_{i}}\otimes_{R} M_{\geq q})_{i+j}.$
It follows that for $j>q$ we get
$$0=\Tor_{i}^{R}(M,K)_{i+j}
=
H_{i}( M \otimes_{R} F_{\lpnt})_{i+j}
=
H_{i}(
M_{\geq q} \otimes_{R} F_{\lpnt})_{i+j} = \Tor_{i}^{R}(M_{\geq q},K)_{i+j}.
$$
Since for $j<q$  we have $(M_{\geq q})_{j}=0$,  we get that
$( M_{\geq q} \otimes_{R} R(-i)^{c_{i}})_{i+j}=0$,
and thus
$\Tor_{i}^{R}(M_{\geq q},K)_{i+j}=0 \text{ for }j<q.$
All in all we proved (iii).

Assume  (iii) holds.
The computation above shows that for integers $j>q$ we have that
$\Tor_{i}^{R}(M,K)_{i+j}=\Tor_{i}^{R}(M_{\geq q},K)_{i+j}=0,$
which shows (ii). This concludes the proof.
\end{proof}

%
%
%
\section{The borderline cases}
\label{borderline}
It is a natural question to characterize the situations where
we have equalities
$\reg^L_R (M)-\reg^L_R (R)=\reg^T_R (M)$
and $\reg^T_R (M)=\reg^L_R (M)+\reg^T_R (K)$ respectively.
Over a standard graded $K$-algebra
the cases that these equalities hold for all finitely generated graded $R$-modules
are easily described.
In fact, Eisenbud and Goto \cite{EIGO}
proved that $\reg^L_R (M)=\reg^T_R (M)$
for all
finitely generated graded $R$-modules $M$ if $R$ is a standard graded polynomial ring.
The next theorem shows that a standard graded polynomial ring
is the only standard graded $K$-algebra with this property.
This results extends also in the module case
an observation in \cite[Corollary 2.8]{JO}.

\begin{thm}
\label{allmodules}
Let $R$ be a standard graded $K$-algebra.
The following statements are equivalent:
\begin{enumerate}
\item
For all finitely generated graded $R$-modules $M$ we have
$\reg^L_R (M)-\reg^L_R (R)=\reg^T_R (M)$;
\item
For all finitely generated graded $R$-modules $M$ we have
$\reg^T_R (M)=\reg^L_R (M)+\reg^T_R (K)$;
\item
For all finitely generated graded $R$-modules $M$ we have
$\reg^T_R (M)=\reg^L_R (M)$;
\item
$R$ is Koszul and $\reg^L_R (R)=0$;
\item
$R=K[x_1,\dots,x_n]$ is a standard graded polynomial ring.
\end{enumerate}
\end{thm}
\begin{proof}
(iv) $\Rightarrow$ (i) , (ii), (iii):
Assume that $R$ is Koszul and $\reg^L_R (R)=0$.
Since $R$ is Koszul we have that $\reg^T_R (K)=0$
by Theorem \ref{koszulavpe}.
Let $M$ be a finitely generated graded $R$-module.
It follows from Theorem \ref{haupt} that
$$
\reg^L_R (M)
=
\reg^L_R (M)-\reg^L_R (R)
\leq
\reg^T_R (M)
\leq
\reg^L_R (M)+\reg^T_R (K)
=
\reg^L_R (M).
$$
Hence $\reg^T_R (M)=\reg^T_L (M)$ in this case.
Thus (i), (ii) and (iii) hold.

(i) $\Rightarrow$ (iv):
Assume that for
all finitely generated graded $R$-modules $M$ we have that
$\reg^L_R (M)-\reg^L_R (R)=\reg^T_R (M)$.
For $M=K$ we get that $\reg^T_R (K)=\reg^L_R (K)-\reg^L_R (R)<\infty$.
It follows from Theorem \ref{koszulavpe} that $\reg^T_R (K)=0$.
Thus $R$ is Koszul and
$\reg^L_R (R)=\reg^L_R (K)-\reg^T_R (K)=0$
where the last equality follows from Remark \ref{goodtoknow}.

(ii) $\Rightarrow$ (iv):
Now we assume that
for all finitely generated graded $R$-modules $M$ we have
$\reg^T_R (M)=\reg^L_R (M)+\reg^T_R (K)$.
For $M=R$ we get that
$0=\reg^T_R (R)=\reg^L_R (R)+\reg^T_R (K)$.
In particular, $\reg^T_R (K)<\infty$.
It follows again from Theorem \ref{koszulavpe} that $\reg^T_R (K)=0$.
Hence $R$ is Koszul and
$\reg^L_R (R)=-\reg^T_R (K)=0$.

(iii) $\Rightarrow$ (iv):
This is shown analogously to the proof of ``(ii) $\Rightarrow$ (iv)''.

(v) $\Rightarrow$ (iv):
If $R=K[x_1,\dots,x_n]$ is a standard graded polynomial ring,
then $R$ is of course Koszul because the Koszul complex
provides a linear free resolution for the $R$-module $K$.
But we also know  $H^i_\mm(R) =0$ for $i \neq n$
and $H^{n}_\mm(R) \cong K[x^{-1}_1,\dots,x^{-1}_{n}](n)$
as $\Z$-graded $R$-modules.
Hence $\reg^L_R (R)=0$.

(iv) $\Rightarrow$ (v):
Next we assume that $R$ is Koszul and $\reg^L_R (R)=0$.
Let $R=S/I$ where
$S=K[x_1,\dots,x_n]$ is a standard graded polynomial ring
and $I \subset S$ is a graded ideal. We also denote by $\mm=(x_1,\dots,x_n)$
the graded maximal ideal of $S$ and
without loss of generality we assume that $I \subseteq \mm^2$ contains no linear forms.
Since the local cohomology of $R$ with respect to $\mm$ as an $R$-module
is isomorphic to the local cohomology of $R$
with respect to $\mm$ as an $S$-module,
we have $\reg^L_S (R)=0$.
For finitely generated graded $S$-modules $M$
we know already that $\reg^L_S (M)=\reg^T_S (M)$
by what we have proved so far. (Use (v) $\Rightarrow$ (iv) and the equivalence of (iii) and (iv).)
Hence $\reg^T_S (R)=0$. But then it follows that $I=(0)$ is the only possibility and thus $R=S$
is a standard graded polynomial ring.
This concludes the proof.
\end{proof}

The latter result shows that
for a standard graded $K$-algebra
the borderline cases of Theorem \ref{haupt}
hold for all finitely generated graded modules only over a polynomial ring.
But it is still a natural question to characterize for an arbitrary positively graded $K$-algebra
which modules have extremal properties with respect to the bounds in
Theorem \ref{haupt}.
Surprisingly we have that for graded modules of finite projective dimension
always the lower inequality of Theorem \ref{haupt} is an equality.

\begin{thm}
\label{finitepd}
Let $R$ be a positively graded $K$-algebra and $M$ be a finitely generated graded $R$-module
such that $\pd_R (M)<\infty$.
Then
$$\reg_{R}^{L}(M) - \reg_{R}^{L} (R) = \reg_{R}^{T}(M).$$
\end{thm}
\begin{proof}
We prove the assertion by induction on $\pd_R (M)$.
Assume first that $\pd_R (M)=0$, then there exist  finitely many $a_{i}\in \Z$ such that
$$
0\to \bigoplus_{i} R(-a_{i}) \to M \to 0
$$
is a  minimal graded free resolution  of $M$ over $R$.
It follows from the definition of $\reg_{R}^{T}$
that
$$
\reg_{R}^{T}(M)=\max\{a_{i}\}.
$$
Moreover, we see that
$$H_{\mm}^{i}(M)_{k-i}=H_{\mm}^{i}(\bigoplus_{j} R(-a_{j}))_{k-i}
=\bigoplus_{j} H_{\mm}^{i}(R)_{k-a_{j}-i},
$$
and thus as desired
$$\reg_{R}^{L}(M)=\reg_{R}^{L} (R) + \reg_{R}^{T}(M).
$$
Assume now $0<\pd_R (M)<\infty$.
Let $F_{0}$ be the first graded free module in the minimal graded free resolution of $M$ over $R$
and let $G_{1}$ be the kernel of the map $F_0 \to M$. Thus we have
the short exact sequence
$$
(*)\qquad 0 \to G_{1} \to F_{0} \to M \to 0.
$$
We have $\pd_R (G_{1})=\pd_R (M) -1$ and hence we can apply the induction hypotheses to $G_1$:
$$
\reg_{R}^{L}(G_{1})
=\reg_{R}^{T}(G_{1})+\reg_{R}^{L}(R).
$$
Because of the definitions of the minimal graded free resolution of $M$ and of $\reg_{R}^{T}$
we have
$$
\reg_{R}^{T}(G_{1})\leq \reg_{R}^{T} (M) + 1.
$$
Let $F_{0}=\bigoplus_{i} R(-a_{i})$ and
$$a:=\max\{a_{i}\}.$$
Thus $\reg_{R}^{L}(F_{0})=a+\reg_{R}^{L}(R)$ and
$a\leq \reg_{R}^{T}(M).$
Now we have to distinguish three cases:

(a)
$\reg_{R}^{T}(G_{1})= \reg_{R}^{T}(M) + 1$:
For $\reg_{R}^{L}(G_{1})$ there exists an integer $j \in \Z$ such that
$H_{\mm}^{j}(G_{1})_{\reg_{R}^{L}(G_{1})-j}\neq 0 $.
It follows from  $(*)$ that
$$\ldots \to H_{\mm}^{j-1}(M)_{\reg_{R}^{L}(G_{1})-1-(j-1)} \to
H_{\mm}^{j}(G_{1})_{\reg_{R}^{L}(G_{1})-j} \to H_{\mm}^{j}(F_{0})_{\reg_{R}^{L}(G_{1})-j}
\to \ldots
$$
Since
$$\reg_{R}^{L}(G_{1})=\reg_{R}^{T}(G_{1})+\reg_{R}^{L}(R)
=\reg_{R}^{T}(M)+1+\reg_{R}^{L}(R) \geq a+1+\reg_{R}^{L}(R)>a+\reg_{R}^{L}(R)$$
we have
$H_{\mm}^{j}(F_{0})_{\reg_{R}^{L}(G_{1})-j}=0$.
Now
$$H_{\mm}^{j-1}(M)_{\reg_{R}^{L}(G_{1})-1-(j-1)}\neq 0$$
because $H_{\mm}^{j-1}(M)_{\reg_{R}^{L}(G_{1})-1-(j-1)}$ maps surjective to
$H_{\mm}^{j}(G_{1})_{\reg_{R}^{L}(G_{1})-j} \neq 0$.
We get
$$\reg_{R}^{L}(M)\geq
\reg_{R}^{L}(G_{1})-1=\reg_{R}^{L}(R)+\reg_{R}^{T}(G_{1})-1=\reg_{R}^{L}(R) +\reg_{R}^{T}(M).$$
By Theorem \ref{haupt} we know already
$$
\reg_{R}^{L}(M)\leq \reg_{R}^{L}(R) +\reg_{R}^{T}(M).
$$
Thus we have equality and the desired assertion follows in this case.

(b)
$\reg_{R}^{T}(G_{1})< \reg_{R}^{T}(M)$:
For the number $a$ as defined as above we have
$$a=\reg_{R}^{T}(M)>\reg_{R}^{T}(G_{1})=\reg_{R}^{L}(G_{1})-\reg_{R}^{L}(R).$$
For the number
$\reg_{R}^{L}(F_{0})$ there exists an $j \in \Z$ such that
$H_{\mm}^{j}(F_{0})_{\reg_{R}^{L}(F_{0})-j}\neq 0 $. By $(*)$
we have the exact sequence
$$\ldots \to H_{\mm}^{j}(G_{1})_{\reg_{R}^{L}(F_{0})-j} \to
H_{\mm}^{j}(F_{0})_{\reg_{R}^{L}(F_{0})-j} \to H_{\mm}^{j}(M)_{\reg_{R}^{L}(F_{0})-j} \to
\ldots\text{ }.$$
Now
$$\reg_{R}^{L}(F_{0})=a+\reg_{R}^{L}(R)=\reg_{R}^{T}(M)+\reg_{R}^{L}(R)>
\reg_{R}^{T}(G_{1})+\reg_{R}^{L}(R)=\reg_{R}^{L}(G_{1}).$$
Thus $H_{\mm}^{j}(G_{1})_{\reg_{R}^{L}(F_{0})-j}=0$
and
$$H_{\mm}^{j}(M)_{\reg_{R}^{L}(F_{0})-j}\neq 0$$
since $H_{\mm}^{j}(F_{0})_{\reg_{R}^{L}(F_{0})-j}$ maps injective  into
$H_{\mm}^{j}(M)_{\reg_{R}^{L}(F_{0})-j}$.
We obtain
$$\reg_{R}^{L}(M) \geq
\reg_{R}^{L}(F_{0})=a+\reg_{R}^{L}(R)=\reg_{R}^{T}(M)+\reg_{R}^{L}(R).$$
It follows again from Theorem \ref{haupt} that
$$\reg_{R}^{L}(M)\leq \reg_{R}^{L}(R) +\reg_{R}^{T}(M).$$
Hence
we have equality and the assertion follows in case (b).

(c) $\reg_{R}^{T}(G_{1})=\reg_{R}^{T}(M)$:
We have for $a$ as defined as above that
$$
a=\reg_{R}^{T}(M)=\reg_{R}^{T}(G_{1}).
$$
For the number $\reg_{R}^{L}(F_{0})$ there exists an integer $j \in \Z$ such that
$H_{\mm}^{j}(F_{0})_{\reg_{R}^{L}(F_{0})-j}\neq 0$.
More precisely, if we write $F_0=R(-a)\oplus F'_0$ for some graded free $R$-module $F'_0$,
then we can assume that
$$
H_{\mm}^{j}(R(-a))_{\reg_{R}^{L}(F_{0})-j}\neq 0
$$
and the induced projection map
$$
\tau_1\colon H_{\mm}^{j}(F_{0})_{\reg_{R}^{L}(F_{0})-j} \to H_{\mm}^{j}(R(-a))_{\reg_{R}^{L}(F_{0})-j}
$$
is surjective.
By $(*)$ we have the exact sequence
$$\ldots \to H_{\mm}^{j}(G_{1})_{\reg_{R}^{L}(F_{0})-j} \to
H_{\mm}^{j}(F_{0})_{\reg_{R}^{L}(F_{0})-j} \to H_{\mm}^{j}(M)_{\reg_{R}^{L}(F_{0})-j} \to
\ldots\text{ }.$$
If
$$\tau_2: H_{\mm}^{j}(G_{1})_{\reg_{R}^{L}(F_{0})-j} \to
H_{\mm}^{j}(F_{0})_{\reg_{R}^{L}(F_{0})-j}$$
would not be surjective,
then $H_{\mm}^{j}(M)_{\reg_{R}^{L}(F_{0})-j} \neq 0$ and it follows that
$$
\reg_{R}^{L}(M) \geq \reg_{R}^{L}(F_{0})=a+\reg_{R}^{L}(R)=\reg_{R}^{T}(M)+\reg_{R}^{L}(R).$$
Again we know from Theorem \ref{haupt} that
$$\reg_{R}^{L}(M)\leq \reg_{R}^{L}(R) +\reg_{R}^{T}(M)$$
and thus we have the desired equality.

It remains to show that indeed $\tau_2$ is not surjective.
Assume for a moment that $\tau_2$ is surjective. Then also the composed map
$$
\tau_3=\tau_1\circ \tau_2 \colon
H_{\mm}^{j}(G_{1})_{\reg_{R}^{L}(F_{0})-j}
\to
H_{\mm}^{j}(R(-a))_{\reg_{R}^{L}(F_{0})-j}
$$
would be surjective. In particular, $\tau_3$ is not the zero map.
Now we consider again
the short exact sequence
$$
0 \to G_{1} \overset{\phi_1}{\to} F_{0} \to M \to 0.
$$
Since $F_0$ was the first module in the minimal graded free resolution of $M$,
we have that $\Image \phi_1 \subseteq \mm F_{0}$.
Since
$\reg_{R}^{T}(G_{1})=\reg_{R}^{T}(M)=a$, we know that
$G_1$ is generated in degrees $\leq a$.
But then for any generator $x$
of $G_1$, it is not possible that $\phi_1(x)$ involves the free generator
corresponding to $R(-a)$ in $F_0$.
In other words, if we compose $\phi_1$ with the natural projection map $\phi_2\colon F_0 \to R(-a)$,
then the induced map $\phi_3 \colon G_1 \to R(-a)$ is the zero map.

Next we observe that the maps
$\tau_1$, $\tau_2$ and $\tau_3$
are induced by $\phi_2, \phi_1$ and $\phi_3$.
Indeed, consider
the  modified Cech-complex $C_{\lpnt}$  (e.g. see \cite[page 130]{BRHE98}).
Then for some graded $R$-modules $W,W'$ and a homogenous map $\psi\colon W \to W'$
we have that
$
H_{\mm}^{j}(W)=H^{i}(W\otimes_{R}C_{\lpnt}),
$
and the natural map $H_{\mm}^{j}(W) \to H_{\mm}^{j}(W')$
corresponds to $\psi\otimes_{R}C_{\lpnt}$.
This implies that the map $\tau_3$ has to be the zero map, because already $\phi_3$
is the zero map. Thus we have a contradiction.
This concludes the proof.
\end{proof}

Now one could hope the converse of Theorem \ref{finitepd} is also true.
But this is not the case as
the next example shows.

\begin{ex}
\label{counterex}
Let $K[x,y]$ be a standard graded polynomial ring in 2 variables
and consider
$$
R=K[x,y]/(x^{2},xy,y^{2}).
$$
Then $R$ is a Koszul algebra since its defining ideal is a monomial ideal generated in degree 2.
(See \cite{FR75}.) $R$ is zero dimensional and thus Cohen--Macaulay.
Let $\omega_{R}$ be the graded canonical module of $R$. Then we have
$$
\reg_{R}^{L} (\omega_{R})-\reg_{R}^{L}(R) = \reg_{R}^{T}(\omega_{R}),
$$
but $\pd_R (\omega_{R})=\infty$.
\end{ex}

\begin{proof}
In the following we identify ideals of $R$ and $S$.
Let $\mm=(x,y)$ be the maximal ideal of $R$.
We  have $\mm=\sqrt{(0)}$  and $\mm^{2}=0$.
For a graded $K$-vector space $W$ we set
$$s(W)=\max\{i \in \Z:W_{i}\neq 0\}.$$
Since $R$ is zero-dimensional and thus also $\omega_{R}$ is zero-dimensional,
it follows from Remark \ref{goodtoknow} that
$
R=H_{\mm}^{0}(R),\ \omega_{R}=H_{\mm}^{0}(\omega_{R}),
$
and for $i>0$ that
$H_{\mm}^{i}(R)=H_{\mm}^{i}(\omega_{R})=0.$
Hence
$$\reg_{R}^{L}(R)=s(R) \text{ and }\reg_{R}^{L}(\omega_{R})=s(\omega_{R}).$$
By the definition of $R$ we have  $s(R)=1$.
By graded local duality we know
$$
\omega_{R}=\Hom_{K}(R,K)
$$
with $(\omega_{R})_i=\Hom_{K}(R_{-i},K)$.
Thus we see that $\dim_K (\omega_{R})_{-1}=2$, $\dim_K (\omega_{R})_{0}=1$ and $\dim_K (\omega_{R})_{i}=0$ for $i\neq -1,0$.
Hence $s(\omega_{R})=0$
and we have
$$
\reg_{R}^{L}(\omega_{R})-\reg_{R}^{L}(R)=-1.
$$
$\omega_{R}$ is a faithful module, thus not all generators of $\omega_{R}$ can be annihilated.
It follows that $\omega_{R}$ is generated in degree $-1$ with $2$ minimal generators.
The minimal graded free resolution of $\omega_{R}$ starts with
$$\ldots\to R^{2}(+1) \to \omega_{R} \to 0.$$
Since $\mm^{2}=0$ in $R$ and the matrices corresponding to the maps in a minimal graded free resolution of $\omega_R$
have entries in $\mm$, we see that $\omega_{R}$ has a $(-1)$-linear resolution.
In particular,
$$\reg_{R}^{T}(\omega_{R})=-1.$$
Thus it follows
$$\reg_{R}^{T}(\omega_{R})=\reg_{R}^{L}(\omega_{R})-\reg_{R}^{L}(R).$$
Assume that $ \pd_R (\omega_{R})<\infty$. Then it follows from the Auslander-Buchsbaum
formula that $\pd_R (\omega_{R})=\depth (R) - \depth (\omega_{R})=0-0=0$.
Hence $\omega_{R}$ would be free which is not possible.
We see that $\pd_R(\omega_{R})=\infty$.
\end{proof}

So it is still interesting to understand better
the modules for which the extremal cases of Theorem \ref{haupt} hold
and we end this section with the following questions.

\begin{quest}
Let $R$ be a positively graded $K$-algebra.
Can one characterize those finitely generated graded $R$-modules $M$
such that
$
\reg_R^L (M)- \reg_R^L (R)=\reg_R^T (M)?
$
\end{quest}

The other inequality
$
\reg_R^L (M)=\reg_R^T (M)+\reg_R^T (K)
$
is only interesting for $R$ a Koszul algebra.
Thus one might ask:
\begin{quest}
Let $R$ be a Koszul-algebra.
Can one characterize those finitely generated graded $R$-modules $M$
such that
$
\reg_R^L (M)=\reg_R^T (M)?
$
\end{quest}

%
%
%
\section{Concluding Examples}
\label{exampleseries}
Let $R$ be a standard graded $K$-algebra and $M$ be a finitely generated graded $R$-module.
In Theorem \ref{haupt} we proved that
$$
\reg_R^L (M) - \reg_R^L (R)
\leq
\reg_R^T (M)
\leq
\reg_R^L (M) + \reg_R^T (K).
$$
We saw that if $\pd_R (M)<\infty$,
then $\reg^L_R (M)-\reg^L_R (R)=\reg^T_R (M)$ is satisfied.
For $R$ Koszul and $M=K$ we see that $\reg^T_R (M)=\reg^L_R (M)$ is true.
Now it is a natural question whether in principle all values between
$\reg^L_R (M)-\reg^L_R (R)$ and $\reg^L_R (M)$ are possible for the number $\reg^T_R (M)$.
We will show that this is true over Koszul algebras.
For this we need at first the following lemma.

\begin{lem}
\label{helper}
Let $R$ be a Koszul-algebra and $M$ be a finitely generated graded $R$-module.
If $M$ has a $(j-1)$-linear resolution, then $\mm M$ has a $j$-linear resolution.
In particular, we have $\reg_{R}^{T} (\mm^{j})=j$ for $j\geq 0$.
\end{lem}
\begin{proof}
We consider the short exact sequence
$$0 \to \mm M \to M \to M/\mm M \to 0$$
and the induced long exact Tor-sequence
$$
\ldots\to \Tor_{i+1}^{R}(M/\mm M,K) \to \Tor_{i}^{R}(\mm M,K) \to \Tor_{i}^{R}(M,K) \to \ldots$$
$$\to \Tor_{1}^{R}(M/\mm M,K) \to \Tor_{0}^{R}(\mm M,K) \to \Tor_{0}^{R}(M,K) \to
\Tor_{0}^{R}(M/\mm M,K) \to 0.$$

Since $M$ has a  $(j-1)$-linear resolution,
we have in particular, that $M$ is generated in degree $j-1$.
The module $M/\mm M$ is a finitely generated graded $K$-vector space.
Hence
$$
M/\mm M\cong \bigoplus K(-j+1)
$$
and this is also an isomorphism of graded $R$-modules.
The minimal graded free resolution of $M/\mm M$ is a direct sum
of the linear minimal graded free resolutions of $K$ shifted by $j-1$.
Thus we see that $M/\mm M$ has an  $(j-1)$-linear resolution.
For $k\neq j-1$ we obtain
$$\Tor_{i}^{R}(M/\mm M,K)_{i+k}=0.$$
Considering again the long exact Tor-sequence above in degree $i+k$
for $k> j$
we get
$$\ldots\to \Tor_{i+1}^{R}(M/\mm M,K)_{i+1+k-1} \to \Tor_{i}^{R}(\mm M,K)_{i+k} \to
\Tor_{i}^{R}(M,K)_{i+k} \to \ldots$$
and therefore
$$\Tor_{i}^{R}(\mm M,K)_{i+k}=0.$$
We have
$\Tor_{i}^{R}(\mm M,K)_{i+k}=0$ for $k<j$ because $\mm M$ is generated in degrees $\geq j$.
Thus we get that
$\mm M$ has a $j$-linear resolution over $R$.

Since $R$ is Koszul, $K=R/\mm$ has a $0$-linear resolution over $R$ which is equivalent to the fact that
$\mm$ has a $1$-linear resolution over $R$.
An induction on $j \geq 1$ yields that $\mm^j$ has a
$j$-linear resolution over $R$.
\end{proof}

\begin{ex}
Let $R$ be a Koszul algebra such that $\depth(R)>0$ and $r=\reg_{R}^{L}(R)>0$.
Then we have for $0<j<r$
that
$$
0=\reg_{R}^{L}(\mm^j)-r < \reg_{R}^{T}(\mm^j)=j <r=\reg_{R}^{L}(\mm^j).
$$
For example consider the $d$-th Veronese subring $S^{(d)}$
of a standard graded polynomial ring $S=K[x_1,\dots,x_n]$
for some integer $d>0$ (i.e. $S^{(d)}$ is the graded $K$-algebra with $(S^{(d)})_i=S_{id}$ for $i\geq 0 $).
For $d\gg 0$ we have that $S^{(d)}$ is Koszul, $\depth(S^{(d)})>0$ and $\reg_{S^{(d)}}^{L} (S^{(d)})>0$.
\end{ex}

\begin{proof}
It follows from Lemma \ref{helper}  that
$$
\reg_{R}^{T} (\mm^{j})=j.
$$
To determine  $\reg_{R}^{L} (\mm^{j})$ we consider the short exact sequence
$$
0 \to \mm^{j} \to R \to R/\mm^{j} \to 0.
$$
The induced long exact local cohomology sequence is
$$0 \to H_{\mm}^{0}(\mm^{j}) \to H_{\mm}^{0}(R) \to H_{\mm}^{0}(R/\mm^{j})\to \ldots$$
$$\to H_{\mm}^{i}(\mm^{j}) \to H_{\mm}^{i}(R) \to H_{\mm}^{i}(R/\mm^{j})\to \ldots.$$
Observe that $R/\mm^{j}$ has finite length, is therefore zero dimensional and we have that
$H_{\mm}^{i}(R/\mm^{j})=0$ for $i>0$.
Since $\depth R>0$ we have
$H_{\mm}^{0}(R)=0$.
Hence
$$
H_{\mm}^{0}(\mm^{j})\subseteq H_{\mm}^{0}(R)=0.
$$
Considering again the long exact local cohomology sequence
we have
$$H_{\mm}^{i}(\mm^{j})\cong H_{\mm}^{i}(R)$$
for $i=0$ and $i>1$.
Moreover, the following sequence is exact:
$$0 \to H_{\mm}^{0}(R/\mm^{j}) \to H_{\mm}^{1}(\mm^{j}) \to H_{\mm}^{1}(R) \to 0.$$
Let $k\geq j-1$. Then we have
$$
H_{\mm}^{0}(R/\mm^{j})_{k}=(R/\mm^{j})_{k}
\begin{cases}
=0 & \text{for $k> j-1$,} \\
\neq 0 & k=j-1.
\end{cases}$$
Thus
$$
H_{\mm}^{1}(R)_{k}=0 \text{ for } k>r
$$
and we see that
$$
\reg_{R}^{L}(\mm^j)=\max\{j, r\}.
$$
For $0<j<r$ we obtain the desired equalities
$$
0=\reg_{R}^{L}(\mm^j)-r < \reg_{R}^{T}(\mm^j)=j <r=\reg_{R}^{L}(\mm^j).
$$
Now let $S=K[x_1,\dots,x_n]$ be a standard graded polynomial ring.
It is well-known that
for $d\gg 0$ the $d$-th veronese subring $S^{(d)}$ of $S$ is Koszul.
(E.g. see \cite{BAFR} or \cite{EIRETO}.)
The number $\reg_{S^{(d)}}^L (S^{(d)})$ coincides
with $\reg_{T}^L (S^{(d)})=\reg_{T}^T (S^{(d)})$ where $T$ is some polynomial ring such that
$S^{(d)}=T/J$ for some graded ideal $J$ containing no linear forms.
But $J$ is generated in degree 2 since $S^{(d)}$ is Koszul. Hence
$\reg_{T}^T (S^{(d)}) \geq 1>0$.
Since $S^{(d)}$ is Cohen--Macaulay of dimension $n$
(e.g. see \cite[Excercise 3.6.21]{BRHE98}) we have in particular $\depth (S^{(d)})>0$.
This shows that we can apply the example to the $K$-algebra $S^{(d)}$.
\end{proof}

\end{document}